\documentclass{article}
\usepackage{gdg-bpblocks}
\usepackage[margin=1.0in]{geometry}

\usepackage{pdfsync}

\author{Andrew Gainer-Dewar \\ Carleton College \\ \url{againerdewar@carleton.edu}
 \and Ira M. Gessel\thanks{Ira M. Gessel's research was partially supported by NSA Grant H98230-10-1-0196 and Simons Foundation Grant \#229238.} \\ Brandeis University \\ \url{gessel@brandeis.edu}}
\title{Enumeration of bipartite graphs and bipartite blocks}

\usepackage{titlesec}
\titleformat{\section}
 {\normalfont\Large\bfseries}{\thesection.}{.8ex}{}
\titleformat{\subsection}
 {\normalfont\large\bfseries}{\thesubsection.}{.8ex}{}
\titleformat{\subsubsection}
 {\normalfont\bfseries}{\thesubsubsection.}{.8ex}{}

\begin{document}
\maketitle

\begin{abstract}
 We use the theory of combinatorial species to count unlabeled bipartite graphs and bipartite blocks (nonseparable or 2-connected graphs). We start with bicolored graphs, which are bipartite graphs that are properly colored in two colors. The two-element group $\mathfrak{S}_{2}$ acts on these graphs by switching the colors, and connected bipartite graphs are orbits of connected bicolored graphs under this action. From first principles we compute the $\mathfrak{S}_{2}$-cycle index for bicolored graphs, an extension of the ordinary cycle index, introduced by Henderson, that incorporates the $\mathfrak{S}_{2}$-action. From this we can compute the
$\mathfrak{S}_{2}$-cycle index for connected bicolored graphs, and then the ordinary cycle index for connected bipartite graphs. The cycle index for connected bipartite graphs allows us, by standard techniques, to  count unlabeled bipartite graphs and unlabeled blocks.

\end{abstract}

\section{Introduction}\label{c:intro}
A \emph{bicolored graph} is a graph of which  each vertex has been assigned one of two colors so that each edge connects vertices of different colors.
A \emph{bipartite graph} is a graph  that admits such a coloring. Given $j$ white and $k$ black vertices, there are $2^{jk}$ ways to join vertices of different colors. Thus the number of (labeled) bicolored graphs on $n$ vertices is
\begin{equation}
\label{eq:bicolor}
b_n = \sum_{i+j=n}\binom ni 2^{ij}.
\end{equation}
 Bipartite graphs are not so easy to count directly. Every connected bicolored graph has exactly two colorings in white and black, so we can count bipartite graphs by relating them to connected bipartite graphs. To do this, we use the \emph{exponential formula} \cite[section 5.1]{stanley:ec2}, which implies that if $f(x)=\sum_{n=1}^\infty f_n x^n/n!$ is the exponential generating function for a class $\mathcal C$ of (labeled) connected graphs then $e^{f(x)}$ is the exponential generating function for graphs all of whose connected components belong to  $\mathcal C$. Conversely, if we know the exponential generating function $g(x)$ for graphs all of whose connected components belong to  $\mathcal C$, then the  exponential generating function for $\mathcal C$ is $\log g(x)$.

 It follows that with $B(x) = \sum_{n=0}^\infty b_n x^n/n!,$
where $b_n$ is given by \cref{eq:bicolor}, the exponential generating function for connected bicolored graphs is $\log B(x)$, the exponential generating function for connected bipartite graphs is $\tfrac12\log B(x)$, and the exponential generating function for bipartite graphs is
$e^{\log B(x)/2}=\sqrt{B(x)}$.

Just as arbitrary graphs may be decomposed into their connected components, arbitrary connected graphs may be decomposed into ``blocks''---maximal $2$-connected (or ``nonseparable'') subgraphs.
Techniques developed by Ford and Uhlenbeck \cite{forduhl:combprob1} were applied by Harary and  Robinson \cite{harrob:bipblocks} to show that the exponential generating function $N(x)$ for labeled $2$-connected bipartite graphs is related to the exponential generating function $P(x) = \frac{1}{2} \log B(x)$ for connected bipartite graphs by the equation $\log P'(x) = N' \pbrac{x P'(x)}$.
This equation suffices to compute the number of labeled bipartite blocks on $n$ vertices and their asymptotics. 

\begin{figure}[htb]
 \centering
 \subfloat[A connected bipartite graph with two distinct bicolorings, one of which is shown]{
   \label{fig:bicoloring:distinct}
   \begin{tikzpicture}
     \node [style=wnode] (1) at (0,0) {};
     \path [draw] (1) -- (0,2) node [style=bnode] (2) {};
     \path [draw] (1) -- (3,2) node [style=bnode] (3) {};
     \path [draw] (1) -- (-3,2) node [style=bnode] (4) {};
   \end{tikzpicture}
 }
 \hspace{1in}
 \subfloat[A connected bipartite graph with just one distinct bicoloring]{
   \label{fig:bicoloring:indistinct}
   \begin{tikzpicture}
     \node [style=wnode] (1) at (0,0) {};
     \path [draw] (1) -- (3,0) node [style=bnode] (2) {};
     \path [draw] (2) -- (3,2) node [style=wnode] (3) {};
     \path [draw] (3) -- (0,2) node [style=bnode] (4) {};
     \path [draw] (4) -- (1);
   \end{tikzpicture}
 }
 \caption{Some connected unlabeled bipartite graphs have two distinct bicolorings, but some have only one}
 \label{fig:bicoloring}
\end{figure}
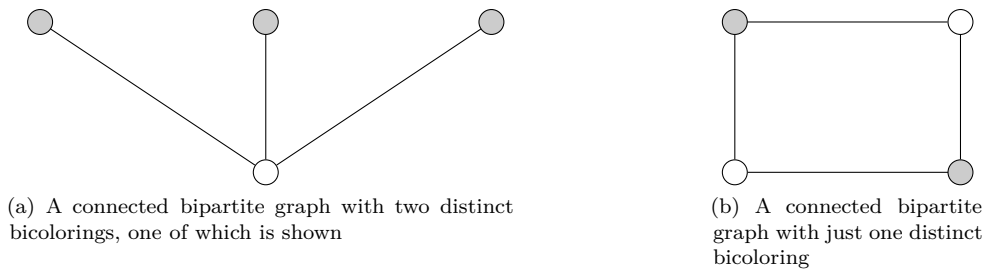

To count unlabeled bipartite graphs we can take a similar approach. It is not too difficult to find the generating function for bicolored graphs from first principles, in a way that is very similar to counting unlabeled graphs (see, e.g., \cite{har:bicolored}).
There is an analogue of the exponential formula for unlabeled graphs (see, for example, \cite[equation~(3.1.1)]{harpalm:graphenum}, \cite[p.~46, equation~(20b) and p.~55, equation~(60\,ii)]{bll:species}, and \cite[p.~29, equation (25) and p.~89, theorem I.5]{fs:anacomb})
so we can easily relate the generating function for all bicolored graphs to that for connected bicolored graphs and the generating function for connected bipartite graphs to that for all bipartite graphs. The difficult step is relating connected bicolored graphs to connected bipartite graphs: some unlabeled connected bipartite graphs can be bicolored in two different ways, and some in only one way, as shown in \cref{fig:bicoloring}.

So instead of simply dividing the number of connected bicolored graphs by 2, as in the labeled case, we must do something more complicated. To deal with this problem, we consider  the two-element group $\symgp{2}$ acting on bicolored graphs by interchanging the colors. We want to count orbits of this group acting on connected bicolored graphs. To do this using Burnside's lemma, we would need to know the number of connected bicolored graphs fixed by each of the elements of $\symgp{2}$. This is not so easy to compute directly; however, it is not hard to compute the corresponding information for the action of $\symgp{2}$ on all bicolored graphs, and by using an extension of the ``unlabeled exponential formula'' we can transfer this information to connected bicolored graphs. Formulas for generating functions for unlabeled bipartite graphs are given in \cref{s:bp}.

In addition to counting unlabeled bipartite graphs, which were counted earlier by Harary and Prins \cite{harprins:bipartite} and by Hanlon \cite{han:bipartite}, we wish to count unlabeled blocks, which have not been previously counted. To accomplish this, we need more than just an enumeration of unlabeled bipartite graphs; we need to compute the cycle index for connected bipartite graphs, a power series in infinitely many variables that generalizes both the labeled and unlabeled enumeration. From the cycle index for connected bipartite graphs, we can use methods developed by Robinson \cite{rob:nonsep} and simplified by Bergeron, Labelle, and Leroux \cite[section 4.2]{bll:species} to count unlabeled bipartite blocks.

To compute the cycle index for connected bipartite graphs we use Henderson's \cite{hend:specfield} extension of Joyal's theory of combinatorial species \cite{bll:species}, which allows us to take account of the $\symgp{2}$-action on bicolored graphs. Our \cref{thm:znbp} gives a formula for the species of bipartite blocks, from which the cycle index for bipartite blocks, and then the ordinary generating function for unlabeled bipartite graphs can be computed.

At the end of their paper on counting labeled bipartite blocks \cite{harrob:bipblocks}, Harary and Robinson stated, ``It is planned to present the counting of unlabeled bipartite blocks in a later communication. Although this is far more difficult than the above labeled enumeration, the cycle index sum methods of \cite{rob:nonsep} can be modified appropriately.'' However, this later communication does not seem to have appeared.

The number of unlabeled bipartite blocks with $n$ vertices for $n\le 24$  is given in \cref{tab:bpblocks}. (Colbourn and Huybrechts~\cite{fggraphs} computed the number of bipartite blocks with at most 14 vertices by generating all connected bipartite graphs and counting those which are 2-connected.)

We would like to thank an anonymous referee for suggestions that improved the presentation of this paper.

\section{The theory of species}\label{c:species}
\subsection{Introduction}\label{s:intro}
Andr\'{e} Joyal \cite{joy:species} introduced the notion of ``species of structures'', which places the  idea of a ``class of labeled objects'' (e.g., trees or permutations) in a categorical setting.
A \emph{species} is a functor from the category $\catname{FinBij}$ of finite sets with bijections to the category $\catname{FinSet}$ of finite sets with set maps\footnote{The use of $\catname{FinSet}$ instead of $\catname{FinBij}$ for the target category is necessary for technical reasons related to quotients.}.
We write $F[A]$ for the image of the set $A$ under the species $F$; for example, if $F$ is the species of bipartite graphs then $F[A]$ is the set of bipartite graphs whose vertices are the elements of the set $A$. We refer the reader to Bergeron, Labelle, and Leroux \cite[\S 1.2]{bll:species} for details; we give here only a brief summary of the facts that we will need from the theory of species.

If $F$ is a species, then for any bijection $\sigma: A\to B$ of finite sets, there is a corresponding bijection $F[\sigma]: F[A]\to F[B]$. Thus if $F[A]$ is the set of graphs with vertex set $A$ and $\sigma$ is a bijection from $A$ to $B$ then for any graph $H\in F[A]$, we obtain $F[\sigma](H)$ by replacing each vertex $a$ of $H$ with $\sigma(a)$. In particular, if $\sigma$ is a bijection from $A$ to itself, then $F[\sigma]$ is a bijection from $F[A]$ to itself, and thus the symmetric group $\symgp{A}$ acts on $F[A]$. The orbits of  $F[A]$ under this action are
``unlabeled'' $F$-structures on~$A$.

We will write $\symgp{n}$ for the symmetric group on the set $\sbrac{n} \defeq {1, 2, \dots, n}$ and for a species $F$ we will write
$F[n]$ for $F[\{1,2,\dots, n\}]$.

Classical enumerative methods frequently use the algebra of generating functions, which record the number of structures of a given size as the coefficients of a formal power series.
To achieve the same goal in species-theoretic analysis, we define an analogous algebraic object which records information related to the action of the permutation groups.
This object is the ``cycle index'' of the species, a symmetric function defined in terms of the power sum symmetric functions $p_i=\sum_j x_j^i$.
(In some accounts of the theory the $p_i$ are taken simply as independent indeterminates.)
\begin{definition}
 \label{def:cycind}
 For a species $F$, we define its \emph{cycle index series} to be the symmetric function
 \begin{equation}
   \label{eq:cycinddef}
   \civars{F}{p_{1}, p_{2}, \dots} = \sum_{n \geq 0} \frac{1}{n!} \biggl( \sum_{\sigma \in \symgp{n}} \fix \pbrac{F \sbrac{\sigma}} p_{\sigma} \biggr),
 \end{equation}
 where $\fix \pbrac{F \sbrac{\sigma}} := \abs{\cbrac{s \in F \sbrac{n} : F \sbrac{\sigma} \pbrac{s} = s}}$, $\sigma_{i}$ is the number of $i$-cycles of $\sigma$, and $p_{\sigma} = p_{1}^{\sigma_{1}} p_{2}^{\sigma_{2}} \dots$.
\end{definition}

It is easy to see%
\footnote{For any group $G$ acting on a set $S$, the number of elements of $S$ fixed by $g\in G$ depends only on the conjugacy class of $g$.}
that $\fix \pbrac{F \sbrac{\sigma}}$ depends only on the cycle type of $\sigma$.
The  cycle types of permutations $\sigma \in \symgp{n}$ are in natural bijective correspondence with integer partitions $\lambda \vdash n$ (that is, weakly decreasing sequences $\pbrac{\lambda_{1}, \lambda_{2}, \dots}$ such that $\sum_{i} \lambda_{i} = n$), and the number of
permutations in $\symgp{n}$ of cycle type $\lambda$ is
$n!/z_{\lambda}$, where if
$\lambda$ has $l_i$ parts equal to $i$ for each $i$, then  $z_\lambda$ is
$1^{l_1}l_1!\, 2^{l_2}l_2!\cdots$.
Thus the contribution to the inner sum in \eqref{eq:cycinddef} from permutations of cycle type $\lambda$ is  $(n!/z_\lambda)\fix (F \sbrac{\lambda}) $,
where
$\fix (F \sbrac{\lambda}) = \fix F (\sbrac{\sigma})$ for any permutation $\sigma$ of cycle type $\lambda$,
and we may write the sum over permutations in \cref{eq:cycinddef} as a sum over partitions:
\begin{equation}
 \label{eq:cycinddefpart}
 \civars{F}{p_{1}, p_{2}, \dots} =
\sum_{n \geq 0} \sum_{\lambda \vdash n} \fix \pbrac{F \sbrac{\lambda}} \frac{p_{\lambda}}{z_{\lambda}},
\end{equation}
where $p_{\lambda} = p_{\lambda_1}p_{\lambda_2}\cdots$.

The cycle index series $\ci{F}$ of the species $F$ captures enough of its structure that we may recover from it both labeled and unlabeled enumerations, though we are concerned here only with unlabeled enumeration, which is given by the following formula:
\begin{theorem}[{\cite[Theorem 8, \S 1.2]{bll:species}}]\label{thm:ciogf}
 The ordinary generating function $\tilde{F} \pbrac{x}$ for unlabeled $F$-structures is given by
 \begin{equation}\label{eq:ciogf}
   \tilde{F} \pbrac{x} = \ci{F} \pbrac[big]{x, x^{2}, x^{3}, \dots}.
 \end{equation}
\end{theorem}

The algebra of cycle indices directly mirrors the combinatorial calculus of species.
Addition, multiplication, and composition of species have natural combinatorial interpretations and correspond directly to addition, multiplication, and plethystic composition of their associated cycle indices.
This last operation is of particular importance:
\begin{definition}
 \label{def:speccomp}
 For two species $F$ and $G$ with $G \sbrac{\varnothing} = \varnothing$, we define their \emph{composition} to be the species $F \circ G$ given by $\pbrac{F \circ G} \sbrac{A} = \prod_{\pi \in P \pbrac{A}} \pbrac{F \sbrac{\pi} \times \prod_{B \in \pi} G \sbrac{B}}$ where $P \pbrac{A}$ is the set of partitions of $A$.
\end{definition}
In other words, the composition $F \circ G$ is the species of $F$-structures of collections of $G$-structures.
\begin{definition}
 \label{def:cipleth}
 Let $f$ and $g$ be cycle indices. Then the \emph{plethysm} $f \circ g$ is the cycle index
 \begin{equation}
   \label{eq:cipleth}
   f \circ g = f \pbrac{g \pbrac{p_{1}, p_{2}, p_{3}, \dots}, g \pbrac{p_{2}, p_{4}, p_{6}, \dots}, \dots},
 \end{equation}
 where $f \pbrac{a, b, \dots}$ denotes the cycle index $f$ with $a$ substituted for $p_{1}$, $b$ substituted for $p_{2}$, and so on.
\end{definition}
This is the same as the definition of plethysm of symmetric functions (see, e.g., Stanley \cite[p.~447]{stanley:ec2}). Plethysm of cycle indices then corresponds exactly to species composition:
\begin{theorem}
 \label{thm:speccompci}
 For species $F$ and $G$ with $G \sbrac{\varnothing} = \varnothing$, the cycle index of their plethysm is
 \begin{equation}
   \label{eq:speccompci}
   \ci{F \circ G} = \ci{F} \circ \ci{G}
 \end{equation}
 where $\circ$ in the right-hand side is as in \cref{eq:cipleth}.
\end{theorem}

Many combinatorial structures admit natural descriptions as compositions of species.
For example, every graph admits a unique decomposition as a (possibly empty) set of (nonempty) connected graphs, so we have the species identity $\specname{G} = \specname{E} \circ \specname{G}^{C}$ where $\specname{E}$ is the species of sets,  $\specname{G}$ the species of graphs,  and $\specname{G}^{C}$ is the species of  connected graphs.

The theory of species may be extended to \emph{virtual species}, which are formal differences of species. All of the operations for species that we have discussed extend in a straightforward way to virtual species. We refer the reader to
\cite[\S 2.5]{bll:species} for details. In particular, two virtual species $F$ and $G$ are \emph{compositional inverses} if $F\circ G=X$ (or equivalently, $G\circ F=X$) where $X$ is the species of singletons, defined by $X[A] = \{A\}$ if $|A|=1$, and $X[A]=\varnothing$ otherwise. We write $F^{\abrac{-1}}$ for the compositional inverse of $F$ if it exists.

\subsection{$\Gamma$-species and quotient species}\label{s:quot}
Burnside's lemma (also known as the Cauchy-Frobenius lemma) is a powerful enumerative tool for counting orbits under a group action. In this section we prove an analogous result for species.

\begin{definition}
 \label{def:gspecies}
 For $\Gamma$ a finite group, a \emph{$\Gamma$-species} $F$ is a combinatorial species $F$ together with an action of $\Gamma$ on $F$-structures which commutes with isomorphisms of those structures.
\end{definition}
For a motivating example, consider the species $k\specname{CG}$ of $k$-colored graphs; the action of $\mathfrak{S}_{k}$ on the colors commutes with relabelings of graphs, so $k\specname{CG}$ is a $\mathfrak{S}_{k}$-species with respect to this action.

From a $\Gamma$-species,  we can construct a  quotient under the action of $\Gamma$:
\begin{definition}
 \label{def:qspecies}
 For $F$ a $\Gamma$-species, define $\nicefrac{F}{\Gamma}$, the \emph{quotient species} of $F$ under the action of $\Gamma$, to be the species of $\Gamma$-orbits of $F$-structures.
\end{definition}

A brief exposition of quotient species may be found in \cite[\S 3.6]{bll:species}, and a more thorough exposition in \cite{bous:species}.

Just as with classical species, we may associate a cycle index to a $\Gamma$-species, following Henderson \cite{hend:specfield}.
\begin{definition}
 \label{def:gcycind}
 For a $\Gamma$-species $F$, we define the $\Gamma$-cycle index $\gci{\Gamma}{F}$: for each $\gamma \in \Gamma$, let
 \begin{equation}
   \gcivars{\Gamma}{F}{\gamma} = \sum_{n \geq 0} \frac{1}{n!} \sum_{\sigma \in \symgp{n}} \fix \pbrac{\gamma \cdot F \sbrac{\sigma}} p_{\sigma} \label{eq:gcycinddef}
 \end{equation}
 with $p_{\sigma}$ as in \cref{eq:cycinddef}.
\end{definition}

We will call such an object (formally a map from $\Gamma$ to the ring $\ringname{Q} \sbrac{\sbrac{p_{1}, p_{2}, \dots}}$ of symmetric functions with rational coefficients in the $p$-basis) a \emph{$\Gamma$-cycle index} even when it is not explicitly the $\Gamma$-cycle index of a $\Gamma$-species.
So the coefficients in the power series count the fixed points of the \emph{combined} action of a permutation and the group element $\gamma$.
Note that, in particular, the classical (``ordinary'') cycle index may be recovered as $\ci{F} = \gcielt{\Gamma}{F}{e}$ for any $\Gamma$-species $F$.

The algebraic relationships between ordinary species and their cycle indices generally extend  to the $\Gamma$-species context. The actions on cycle indices of $\Gamma$-species addition and multiplication are exactly as in the ordinary species case considered $\Gamma$-componentwise.
The action of composition, which in ordinary species corresponds to plethysm of cycle indices, can also be extended:
\begin{definition}
 \label{def:gspeccomp}
 For two $\Gamma$-species $F$ and $G$, define their \emph{composition} to be the $\Gamma$-species $F \circ G$ with structures given by $\pbrac{F \circ G} \sbrac{A} = \prod_{\pi \in P \pbrac{A}} \pbrac{F \sbrac{\pi} \times \prod_{B \in \pi} G \sbrac{B}}$ where $P \pbrac{A}$ is the set of partitions of $A$ and where $\gamma \in \Gamma$ acts on a $\pbrac{F \circ G}$-structure by acting on the $F$-structure and the $G$-structures independently.
\end{definition}
A formula similar to that of \cref{thm:speccompci} requires a definition of the plethysm of $\Gamma$-symmetric functions, here taken from Henderson \cite[\S 3]{hend:specfield}.
\begin{definition}
 \label{def:gcipleth}
 For two $\Gamma$-cycle indices $f$ and $g$, their \emph{plethysm} $f \circ g$ is a $\Gamma$-cycle index defined by
 \begin{equation}
   \pbrac{f \circ g} \pbrac{\gamma} = f \pbrac{\gamma} \pbrac{g \pbrac{\gamma} \pbrac{p_{1}, p_{2}, p_{3}, \dots}, g \pbrac[big]{\gamma^{2}} \pbrac{p_{2}, p_{4}, p_{6}, \dots}, \dots}.
   \label{eq:gcipleth}
 \end{equation}
\end{definition}
This definition of $\Gamma$-cycle index plethysm is then indeed the correct operation to pair with the composition of $\Gamma$-species:
\begin{theorem}[{\cite[Theorem 3.1]{hend:specfield}}]
 \label{thm:gspeccompci}
 If $A$ and $B$ are $\Gamma$-species and $B \pbrac{\varnothing} = \varnothing$, then
 \begin{equation}
   \label{eq:gspeccompci}
   \gci{\Gamma}{A \circ B} = \gci{\Gamma}{A} \circ \gci{\Gamma}{B}.
 \end{equation}
\end{theorem}

Recall from \cref{eq:cycinddef} that, to compute the cycle index of a species, we need to enumerate the fixed points of each $\sigma \in \symgp{n}$.
To count fixed points in the quotient species $\nicefrac{F}{\Gamma}$ we need to count the fixed $\Gamma$-orbits of $\sigma$ in $F$ under commuting actions of $\symgp{n}$ and $\Gamma$ (that is, under an $\pbrac{\symgp{n} \times \Gamma}$-action).
This may be accomplished by the following generalization of Burnside's lemma  \cite{rob:duality}. (A  more general result appears in \cite[Theorem 4.2b]{deb:polya}.)
\begin{lemma}
 \label{thm:robcount}
 If $\,\Gamma$ and $\Delta$ are finite groups and $S$ is a set with a $\pbrac{\Gamma \times \Delta}$-action, then for any $\delta \in \Delta$ the number of $\Gamma$-orbits fixed by $\delta$ is $\frac{1}{\abs{\Gamma}} \sum_{\gamma \in \Gamma} \fix \pbrac{\gamma, \delta}$.
\end{lemma}
Applying
\cref{thm:robcount} to \cref{def:gcycind} yields a formula for the cycle index of a quotient species in terms of the
$\Gamma$-cycle index. An equivalent result was given by Bousquet \cite[\S 2.2.3]{bous:species}.
\begin{theorem}\label{thm:qsci}
 For a $\Gamma$-species $F$, the ordinary cycle index of the quotient species $\nicefrac{F}{\Gamma}$ is given by
 \begin{equation}
   \label{eq:quotcycind}
   \ci{F / \Gamma} =  \frac{1}{\abs{\Gamma}} \sum_{\gamma \in \Gamma} \gcielt{\Gamma}{F}{\gamma}.
     \end{equation}
\end{theorem}

We will use the notation  $\qgci{\Gamma}{F}$ for  $\frac{1}{\abs{\Gamma}} \sum_{\gamma \in \Gamma} \gcielt{\Gamma}{F}{\gamma}$.

When $\Gamma$ is a symmetric group $\symgp{n}$, as in our applications, we may represent the $\Gamma$-cycle index as a symmetric function in two sets of variables which is homogeneous of degree $n$ in the second set of variables; with this approach, $\Gamma$-cycle index plethysm corresponds to the operation of ``inner plethysm in $y$'' studied by Travis \cite{travis:inpleth}.

\section{The species of bipartite blocks}\label{c:bpblocks}
\subsection{Introduction}\label{s:bpintro}

\begin{definition}
 \label{def:bcgraph}
 A \emph{bicolored graph} is a graph, each vertex of which has been assigned one of two colors (here, black and white) such that each edge connects vertices of different colors.
 A \emph{bipartite graph} (sometimes called \emph{bicolorable}) is a graph which admits such a coloring.
\end{definition}

There is an extensive literature about bicolored and bipartite graphs, including enumerative results for bicolored graphs \cite{har:bicolored}, bipartite graphs both allowing \cite{han:bipartite} and prohibiting \cite{harprins:bipartite} isolated points, and bipartite blocks \cite{harrob:bipblocks}.
However, the enumeration of bipartite blocks has been accomplished previously only in the labeled case.
By considering the problem in light of the theory of $\Gamma$-species, we develop a more systematic understanding of the structural relationships between these various classes of graphs, which allows us, in particular, to enumerate all of them in their unlabeled forms.

Throughout this chapter, we denote by $\specname{BC}$ the species of bicolored graphs and by $\specname{BP}$ the species of bipartite graphs.
The prefix $\specname{C}$ will indicate the connected analogue of such a species, so $\specname{CBP}$ is the species of connected bipartite graphs.

We are motivated by the graph-theoretic fact that each \emph{connected} bipartite graph has exactly two bicolorings, and may be identified with an orbit of connected bicolored graphs under the action of $\symgp{2}$ where the nontrivial element $\tau$ reverses all vertex colors.
We will hereafter treat all the various species of bicolored graphs as $\symgp{2}$-species with respect to this action and use the theory developed in \cref{s:quot} to pass to bipartite graphs.

\subsection{Bicolored graphs}\label{s:bcgraph}
We begin our investigation by directly computing the $\symgp{2}$-cycle index for the species $\specname{BC}$ of bicolored graphs with the color-reversing $\symgp{2}$-action described previously.
We will then use various methods from the species algebra of \cref{c:species} to pass to  other related species. To compute the $\symgp{2}$-cycle index $\gci{\symgp{2}}{\specname{BC}}$ we compute separately $\gcielt{\symgp{2}}{\specname{BC}}{e}$ and $\gcielt{\symgp{2}}{\specname{BC}}{\tau}$.

\subsubsection{Computing $\gcielt{\symgp{2}}{\specname{BC}}{e}$}\label{ss:ecibc}
For each $n > 0$ and each permutation $\pi \in \symgp{n}$, we must count bicolored graphs on $\sbrac{n}$ for which $\pi$ is a color-preserving automorphism.
To simplify some future calculations, we omit empty graphs and define $\specname{BC} \sbrac{\varnothing} = \varnothing$.
We note that the \emph{number} of such graphs in fact depends only on the cycle type $\lambda \vdash n$ of the permutation $\pi$, so we can use the cycle index formula in \cref{eq:cycinddefpart} interpreted as a $\Gamma$-cycle index identity.

Fix some $n \geq 0$ and let $\lambda \vdash n$.
We wish to count bicolored graphs for which a chosen permutation $\pi$ of cycle type $\lambda$ is a color-preserving automorphism.
Each cycle of the permutation must correspond to a monochromatic subset of the vertices, so we may construct graphs by drawing bicolored edges into a given colored vertex set.
If we draw some particular bicolored edge, we must also draw every other edge in its orbit under $\pi$ if $\pi$ is to be an automorphism of the graph.
Moreover, every bicolored graph for which $\pi$ is an automorphism may be constructed in this way.
Therefore, we direct our attention first to counting these edge orbits for a fixed coloring; we will then count colorings with respect to these results to get our total cycle index.

\begin{figure}[htb]
 \centering

 \begin{tikzpicture}
   \GraphInit[vstyle=Hasse]

   \begin{scope}[xshift=-3cm,rotate=30]
     \SetUpEdge[style=cycedge]
     \grCycle[RA=2,prefix=a]{4}
   \end{scope}

   \begin{scope}[xshift=+3cm,rotate=90]
     \SetUpEdge[style=cycedge]
     \grCycle[RA=2,prefix=b]{2}
     \AddVertexColor{black!20}{b0,b1}
   \end{scope}

   \SetUpEdge[style=dashed]
   \EdgeDoubleMod{a}{4}{0}{1}{b}{2}{0}{1}{4}

   \SetUpEdge[style=solid]
   \Edge[label={$e$}](a1)(b1)
 \end{tikzpicture}
 \caption[Example edge-orbit of a color-preserving automorphism]{An edge $e$ (solid) between two cycles of lengths $4$ and $2$ in a permutation and that edge's orbit (dashed)}
 \label{fig:exbcecycle}
\end{figure}
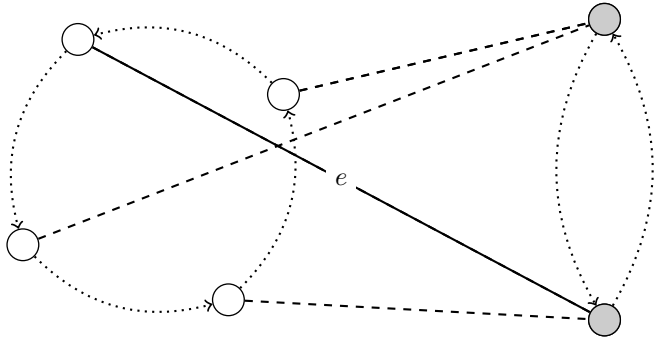

Consider an edge connecting two cycles of lengths $m$ and $n$; the length of its orbit under the permutation is $\lcm \pbrac{m, n}$, so the number of such orbits of edges between these two cycles is $mn / \lcm \pbrac{m, n} = \gcd \pbrac{m, n}$.
For an example in the case $m = 4, n = 2$, see \cref{fig:exbcecycle}.
The number of orbits for a fixed coloring is then $\sum \gcd \pbrac{m, n}$ where the sum is over the multisets of all cycle lengths $m$ of white cycles and $n$ of black cycles in the permutation $\pi$.
We may then construct any possible graph fixed by our permutation by making a choice of a subset of these cycles to fill with edges, so the total number of such graphs is $\prod 2^{\gcd \pbrac{m, n}}$ for a fixed coloring.

We now turn our attention to the possible colorings of the graph which are compatible with a permutation of specified cycle type $\lambda$.
We split our partition into two subpartitions, writing $\lambda = \mu \cup \nu$, where partitions are treated as multisets and $\cup$ is the multiset union, and  $\mu$  corresponds to the white cycles and $\nu$ the black.
Then the total number of graphs fixed by such a permutation with a specified decomposition is
\begin{equation}
 \label{eq:efixedbcgraphs}
 \fix \pbrac{\mu, \nu} = \prod_{\substack{i \in \mu \\ j \in \nu}} 2^{\gcd \pbrac{i, j}}
\end{equation}
where the product is over the elements of $\mu$ and $\lambda$ taken as multisets.
Suppose that the multiplicities of the part $i$ in the partitions $\lambda$, $\mu$, and $\nu$ are $l_i$, $m_i$, and $n_i$, respectively.
Then the $l_i$ $i$-cycles of a permutation of cycle type $\lambda$ can be colored so that  $m_i$ are white and $n_i$ are black in  $l_i!/\pbrac{m_i!\, n_i!}$ ways.
So in all there are $\prod_{i} l_{i}! / \pbrac{m_{i}!\, n_{i}!} = z_{\lambda} / \pbrac{z_{\mu} z_{\nu}}$ colorings associated with $\mu$ and $\nu$, and
\begin{equation*}
 \fix \pbrac{\lambda} = \frac{z_{\lambda}}{z_{\mu} z_{\nu}} \fix \pbrac{\mu, \nu} = \sum_{\mu \cup \nu = \lambda} \frac{z_{\lambda}}{z_{\mu} z_{\nu}} \prod_{\substack{i \in \mu \\ j \in \nu}} 2^{\gcd \pbrac{i, j}}.
\end{equation*}

Thus we obtain a formula for $\gcielt{\symgp{2}}{\specname{BC}}{e}$.
\begin{theorem}
 Let $\specname{BC}$ denote the $\symgp{2}$-species of bicolored graphs with the color-switching action of $\symgp{2}$.
 Then the element of the $\symgp{2}$-cycle index $\gci{\symgp{2}}{\specname{BC}}$ associated to $e$ is given by
 \begin{equation}
   \label{eq:ecibc}
   \gcielt{\symgp{2}}{\specname{BC}}{e} = \sum_{n > 0} \sum_{\substack{\mu, \nu \\ \mu \cup \nu \vdash n}} \frac{p_{\mu \cup \nu}}{z_{\mu} z_{\nu}} \prod_{i, j} 2^{\gcd \pbrac{\mu_{i}, \nu_{j}}}.
 \end{equation}
\end{theorem}

Explicit formulas for the generating function for unlabeled bicolored graphs were obtained by Harary \cite{har:bicolored} using conventional P\'{o}lya-theoretic methods.
Conceptually, our enumeration largely mirrors his.
Harary uses the classical cycle index of the line group\footnote{The \emph{line group} of a graph is the group of permutations of edges induced by permutations of vertices.} of the complete bicolored graph of which any given bicolored graph is a spanning subgraph.
He then enumerates orbits of edges under these groups using the P\'{o}lya enumeration theorem.

\subsubsection{Calculating $\gcielt{\symgp{2}}{\specname{BC}}{\tau}$}\label{ss:tcibc}
Recall that the nontrivial element of $\tau \in \symgp{2}$ acts on bicolored graphs by reversing all colors.

We again consider the cycles in the vertex set $\sbrac{n}$ induced by a permutation $\pi \in \symgp{n}$ and use the partition $\lambda$ corresponding to the cycle type of $\pi$ for bookkeeping.
We then wish to count bicolored graphs on $\sbrac{n}$ for which $\tau \cdot \pi$ is an automorphism, which is to say that $\pi$ itself is a color-\emph{reversing} automorphism.
The number of bicolored graphs for which $\pi$ is a color-reversing automorphism depends only on the cycle type $\lambda$.
Each cycle of vertices must be color-alternating and hence of even length, so the partition $\lambda$ must have only even parts.
Once this condition is satisfied, edges may be drawn either within a single cycle or between two cycles, and as before if we draw in any edge we must draw in its entire orbit under $\pi$ (since $\pi$ is to be an automorphism of the underlying graph).
Moreover, all graphs for which $\pi$ is a color-reversing automorphism with a fixed coloring may be constructed in this way, so it suffices to count such edge orbits and then consider how colorings may be assigned.

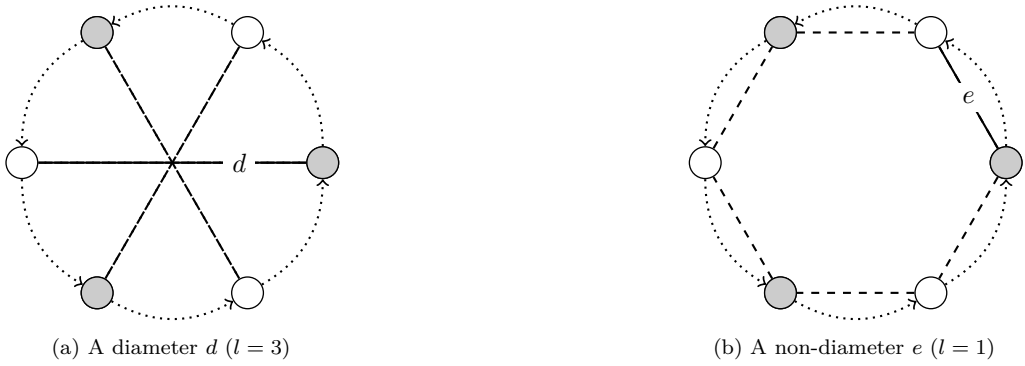
\begin{figure}[htb]
 \centering
 \subfloat[A diameter $d$ ($l = 3$)]{\makebox[.45\textwidth]{
     \label{fig:exbctincycd}
     \begin{tikzpicture}
       \GraphInit[vstyle=Hasse]

       \SetUpEdge[style=cycedge]
       \grCycle[RA=2,prefix=a]{6}
       \AddVertexColor{black!20}{a0,a2,a4}

       \SetUpEdge[style=dashed]
       \EdgeInGraphMod{a}{6}{3}{0}

       \SetUpEdge[style=solid]
       \Edge[label={$d$},style={pos=.25}](a0)(a3)
     \end{tikzpicture}
   }}
 \hfill
 \subfloat[A non-diameter $e$ ($l = 1)$]{\makebox[.45\textwidth]{
     \label{fig:exbctincyce}
     \begin{tikzpicture}
       \GraphInit[vstyle=Hasse]
       \SetUpEdge[style=cycedge]
       \grCycle[RA=2,prefix=b]{6}
       \AddVertexColor{black!20}{b0,b2,b4}

       \SetUpEdge[style=dashed]
       \EdgeInGraphMod{b}{6}{1}{0}

       \SetUpEdge[style=solid]
       \Edge[label={$e$}](b0)(b1)
     \end{tikzpicture}
   }}
 \caption[Two example edge-orbits in a color-reversing automorphism]{Both types of intra-cycle edges and their orbits on a typical color-alternating $6$-cycle}
 \label{fig:exbctincyc}
\end{figure}

We first determine the number of orbits of edges within a cycle of length $2n$; we hereafter describe such a cycle as having \emph{semilength} $n$.
There are exactly $n^{2}$ possible white-black edges in such a cycle. If $n$ is even, then every edge lies in an orbit of size $2n$, so there are $n^2\!/(2n) = n/2$ orbits of edges.
If $n$ is odd, there are $n$ edges joining diametrically opposed vertices, which have oppositive colors. These $n$  edges are all in the same orbit. (See \cref{fig:exbctincycd} for an illustration of these edges.)
The remaining $n^2-n$ edges are in orbits of size $2n$, so there are $(n^2-n)/(2n) = (n-1)/2$ of these orbits. (See \cref{fig:exbctincyce} for an illustration of these edges.)
Thus the total number of orbits for $n$ odd is $(n+1)/2$. In either case, the number of orbits is  $\ceil{n/2}$.

Now consider an edge drawn between two cycles of semilengths $m$ and $n$.
The total number of possible white-black edges is $2mn$, each of which has an orbit length of $\lcm \pbrac{2m, 2n} = 2 \lcm \pbrac{m, n}$.
Hence, the total number of orbits is $2mn / \pbrac{2 \lcm \pbrac{m, n}} = \gcd \pbrac{m, n}$.

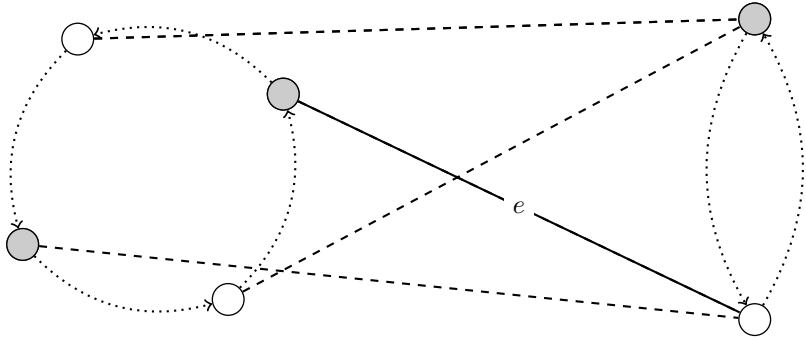
\begin{figure}[htb]
 \centering
 \begin{tikzpicture}
   \GraphInit[vstyle=Hasse]

   \begin{scope}[xshift=-4cm,rotate=30]
     \SetUpEdge[style=cycedge]
     \grCycle[RA=2,prefix=a]{4}
     \AddVertexColor{black!20}{a0,a2}
   \end{scope}

   \begin{scope}[xshift=4cm,rotate=90]
     \SetUpEdge[style=cycedge]
     \grCycle[RA=2,prefix=b]{2}
     \AddVertexColor{black!20}{b0}
   \end{scope}

   \SetUpEdge[style=dashed]
   \EdgeDoubleMod{a}{4}{1}{1}{b}{2}{0}{1}{4}

   \SetUpEdge[style=solid]
   \Edge[label={$e$}](a0)(b1)
 \end{tikzpicture}
 \caption[Another example edge-orbit of a color-reversing automorphism]{An edge $e$ and its orbit between color-alternating cycles of semilengths $2$ and $1$}
 \label{fig:exbctbtwcyc}
\end{figure}

All together, then, the number of orbits for a fixed coloring of a permutation of cycle type $2 \lambda$ (the partition obtained by doubling every part of $\lambda$) is $\sum_{i} \ceil{\lambda_{i}/{2} }+ \sum_{i < j} \gcd \pbrac{\lambda_{i}, \lambda_{j}}$.
All valid bicolored graphs for a fixed coloring for which $\pi$ is a color-preserving automorphism may be obtained uniquely by making some choice of a subset of this collection of orbits, just as in \cref{ss:ecibc}.
Thus, the total number of possible graphs for a given vertex coloring is
\begin{equation}
 \label{eq:tfixedbcgraphs}
 \prod_{i} 2^{\ceil{\lambda_i/2}} \prod_{i < j} 2^{\gcd \pbrac{\lambda_{i}, \lambda_{j}}},
\end{equation}
independent of the choice of coloring.
For a partition $2\lambda$ with $l \pbrac{\lambda}$ cycles, there are then $2^{l \pbrac{\lambda}}$ colorings compatible with our requirement that each cycle is color-alternating, which we multiply by \eqref{eq:tfixedbcgraphs} to obtain the total number of graphs for all permutations $\pi$ with cycle type $2 \lambda$.

Thus, we obtain a formula for $\gcielt{\symgp{2}}{\specname{BC}}{\tau}$:
\begin{theorem}
 Let $\specname{BC}$ denote the $\symgp{2}$-species of bicolored graphs with the color-switching action of $\symgp{2}$.
 Then the element of the $\symgp{2}$-cycle index $\gci{\symgp{2}}{\specname{BC}}$ associated to $\tau$ is given by
 \begin{equation}
   \label{eq:tcibc}
   \gcielt{\symgp{2}}{\specname{BC}}{\tau} = \sum_{\substack{n > 0 \\ \text{\rm $n$ even}}} \sum_{\lambda \vdash n/2} 2^{l \pbrac{\lambda}} \frac{p_{2 \lambda}}{z_{2 \lambda}} \prod_{i} 2^{\ceil{\lambda_i/2}} \prod_{i < j} 2^{\gcd \pbrac{\lambda_{i}, \lambda_{j}}}.
 \end{equation}
\end{theorem}

\subsection{Connected bicolored graphs}\label{s:cbc}
As noted in the introduction to this section, to pass from bicolored graphs to bipartite graphs by taking a quotient under the color-reversing action of $\symgp{2}$, we must work in the connected case.
Thus, we must first pass from the $\symgp{2}$-species $\specname{BC}$ of bicolored graphs to the $\symgp{2}$-species $\specname{CBC}$ of connected bicolored graphs.

Every graph may be decomposed uniquely into (and thus species-theoretically identified with) the set of its connected components.
Thus, at the species level, we have that
\begin{equation} \specname{BC} = \specname{E}^{+} \circ \specname{CBC}, \label{eq:bcdecomp} \end{equation}
where $\specname{BC}$ is the species of nonempty bicolored graphs, $\specname{CBC}$ is the species of nonempty connected bicolored graphs, and $\specname{E}^{+} = \specname{E} - 1$ is the species of nonempty sets.

Reversing the colors of a bicolored graph is done simply by reversing the colors of each of its connected components independently; this action has no effect on the structure of the collection of these components.
We may extend $\specname{E}^{+}$ to an $\symgp{2}$-species by applying the trivial action; then \cref{eq:bcdecomp} may be interpreted as an isomorphism of $\symgp{2}$-species.

To use the decomposition in \cref{eq:bcdecomp} to derive the $\symgp{2}$-cycle index for $\specname{CBC}$, we must invert the $\symgp{2}$-species composition into $\specname{E}^{+}$.
We write $\con := \pbrac{\specname{E}^+}^{\abrac{-1}}$ to denote the virtual species that is the inverse of $\specname{E}^+$ with respect to composition of species, following the notation of \cite{labelle:pssing}; we also let $\con$ denote the virtual $\symgp{2}$-species which is the compositional inverse of the $\symgp{2}$-species $\specname{E}^{+}$ with the trivial action.

We can derive from \cite[\S 2.5, equation~(58c)]{bll:species} a formula for the cycle index $\ci{\con}$ of this virtual species.
It is then straightforward to show (by consideration of \cref{eq:gcipleth}) that each term of $\gci{\symgp{2}}{\con}$ is equal to $\ci{\con}$.
Thus, we have that
\begin{equation}
 \label{eq:clogci}
 \gcielt{\symgp{2}}{\con}{\gamma}= \sum_{k \geq 1} \frac{\mu \pbrac{k}}{k} \log \pbrac{1 + p_{k}}
\end{equation}
for each $\gamma \in \symgp{2}$, where $\mu$ is the integer M\"{o}bius function.
(This cycle index series is sometimes known as the ``combinatorial logarithm''.)

We can then rewrite \cref{eq:bcdecomp} as
\begin{equation}
 \label{eq:bcdecompinv}
 \specname{CBC} = \con \circ \specname{BC}.
\end{equation}
Translating \cref{eq:bcdecompinv} then gives us a formula for the $\symgp{2}$-cycle index of $\specname{CBC}$.

\begin{theorem}
 Let $\specname{BC}$ denote the $\symgp{2}$-species of bicolored graphs and $\specname{CBC}$ the $\symgp{2}$-species of connected bicolored graphs, both with the color-switching action of $\symgp{2}$.
 Additionally, let $\gci{\symgp{2}}{\con}$ denote the combinatorial logarithm series given in \cref{eq:clogci}.
 Then the $\symgp{2}$-cycle indices of $\specname{BC}$ and $\specname{CBC}$ are related by
 \begin{equation} \gci{\symgp{2}}{\specname{CBC}} = \gci{\symgp{2}}{\con} \circ \gci{\symgp{2}}{\specname{BC}}. \label{eq:zcbcdecomp} \end{equation}
\end{theorem}

Note that we could have avoided the use of virtual species by performing the inversion at the level of cycle indices.

\subsection{Bipartite graphs}\label{s:bp}
As we previously observed, connected bipartite graphs are naturally identified with orbits of connected bicolored graphs under the color-reversing action of $\symgp{2}$.
Thus,
\begin{equation*}
 \specname{CBP} = \faktor{\specname{CBC}}{\symgp{2}}.
\end{equation*}
By application of \cref{thm:qsci}, we can then directly compute the cycle index of $\specname{CBP}$ in terms of previous results.

\begin{theorem}
 Let $\specname{CBP}$ denote the species of connected bipartite graphs and $\specname{CBC}$ the $\symgp{2}$-species of connected bicolored graphs.
 Their cycle indices are related by
 \label{thm:zcbp}
 \begin{equation}
   \ci{\specname{CBP}} = \qgci{\symgp{2}}{\specname{CBC}} = \frac{1}{2} \pbrac{\gcielt{\symgp{2}}{\specname{CBC}}{e} + \gcielt{\symgp{2}}{\specname{CBC}}{\tau}}.
 \end{equation}
\end{theorem}

Since a bipartite graph is a set of connected bipartite graphs, we have $\specname{BP} = \specname{E} \circ \specname{CBP}$, and this gives the formula for the cycle index for bipartite graphs.

\begin{theorem}
 Let $\specname{BP}$ denote the species of bipartite graphs, $\specname{CBP}$ the species of connected bipartite graphs, and $\specname{E}$ the species of sets.
 Their cycle indices are related by
 \label{thm:zbp}
 \begin{equation}
   \ci{\specname{BP}} = \ci{\specname{E}} \circ \ci{\specname{CBP}}.
 \end{equation}
\end{theorem}

\cref{thm:zbp} allows us to compute the number of unlabeled bipartite graphs with $n$ vertices.
However,
we can find a computationally more efficient formula for bipartite graphs using only ordinary generating functions, rather than cycle indices.
Specifically, let
\begin{align*}
f_e(x)&= 1+\gcielt{\symgp{2}}{\specname{BC}}{e}(x,x^2,x^3,\dots)\\
f_\tau(x)&=1+ \gcielt{\symgp{2}}{\specname{BC}}{\tau}(x,x^2,x^3,\dots)\\
g_e(x) &=\gcielt{\symgp{2}}{\specname{CBC}}{e}(x,x^2,x^3,\dots)\\
g_\tau(x)&=\gcielt{\symgp{2}}{\specname{CBC}}{\tau}(x,x^2,x^3,\dots)\\
c(x)&= \tilde{Z}_{\specname{CBP}}(x) = \ci{\specname{CBP}}(x,x^2, x^3,\dots)\\
b(x)&=\tilde{Z}_{\specname{BP}}(x) = \ci{\specname{BP}}(x,x^2, x^3,\dots).
\end{align*}
Then $c(x)$ is the ordinary generating function for connected bipartite graphs and $b(x)$ is the ordinary generating function for bipartite graphs. We have formulas for $f_e(x)$ and $f_\tau(x)$ as   sums over partitions,
\begin{align*}
f_e(x) &= \sum_{n=0}^\infty x^n \sum_{\substack{\mu, \nu \\ \mu \cup \nu \vdash n}}
   \frac{1}{z_{\mu} z_{\nu}} \prod_{i, j} 2^{\gcd \pbrac{\mu_{i}, \nu_{j}}}\\
 f_\tau(x) &= \sum_{\text{\rm $n$ even}} x^n\sum_{\lambda \vdash n/2}
  \frac{2^{l \pbrac{\lambda}}}{z_{2 \lambda}} \prod_{i} 2^{\ceil{\lambda_i/2}} \prod_{i < j} 2^{\gcd \pbrac{\lambda_{i}, \lambda_{j}}},
\end{align*}
and $g_e(x)$ and $g_\tau(x)$ are related to $f_e(x)$ and $f_\tau(x)$ by
\begin{align*}
f_e(x)&=\exp\biggl(\sum_{k=1}^\infty \frac{g_e(x^k)}{k}\biggr)\\
f_\tau(x)&=\exp\biggl( \sum_{k=0}^\infty \frac{g_\tau(x^{2k+1})}{2k+1}
 +\sum_{k=1}^\infty \frac{g_e(x^{2k})}{2k}
  \biggr),
\end{align*}
which may be inverted to give
\begin{align}
g_e(x)&=\sum_{k=1}^\infty \frac{\mu(k)}{k} \log f_e(x^k)\notag\\
g_\tau(x)&= \sum_{k=0}^\infty
 \frac{\mu(2k+1)}{2k+1} \log f_\tau(x^{2k+1})
+\sum_{k=1}^\infty \frac{\mu(2k)}{2k} \log f_e(x^{2k}),\label{eq:nohanlon}
\end{align}
where $\mu$ is the M\"obius function.
Finally, $c(x)=\tfrac12(g_e(x) + g_\tau(x))$ and $b(x) = \exp\bigl(\sum_{k=1}^\infty c(x^k)/k\bigr)$.

These calculations are essentially the same as  Hanlon's \cite{han:bipartite}, though he does not have our \cref{eq:nohanlon}.
Unlabeled bipartite graphs were first counted, using a different approach, by Harary and Prins \cite{harprins:bipartite}.

In order to count bipartite blocks, which we accomplish in the next section, we do need the entire cycle index.

\subsection{Nonseparable graphs}\label{s:nbp}
We now turn our attention to the notions of block decomposition and nonseparable graphs.
A graph is said to be \emph{nonseparable} if it is vertex-$2$-connected (that is, if there exists no vertex whose removal disconnects the graph); every connected graph then has a canonical ``decomposition''\footnote{Note that this decomposition does not actually partition the vertices, since many blocks may share a single cut-point.
}
into maximal nonseparable subgraphs, often shortened to \emph{blocks}.
In the spirit of our previous notation, we we will denote by $\specname{NBP}$ the species of nonseparable bipartite graphs, our object of study.

The basic principles of block enumeration in terms of automorphisms and cycle indices of permutation groups were first identified and exploited  by Robinson \cite{rob:nonsep}.
In \cite[\S 4.2]{bll:species}, a theory relating a  species $B$ of nonseparable graphs to the species $C_{B}$ of connected graphs whose blocks are in $B$ is developed using similar principles.

We extract two particular results, appearing as \cite[equations~4.2.27 and 4.2.26a]{bll:species}. We note that the \emph{derivative} $F'$ of a species $F$ \cite[pp.~47--49]{bll:species} is defined by
$F'[A]= F[A\cup \{*\}]$, where $*$ is not in $A$, and its cycle index satisfies
$Z_{F'} = {\partial Z_F}/{\partial p_1}$. The \emph{pointing} $F^{\bullet}$ of $F$ \cite[\S 2.1]{bll:species} is $XF'$. Thus an $F^{\bullet}$-structure on the set $A$ is an element of $F[A]$ together with a distinguished element of $A$.

\begin{theorem}
 \label{thm:blocks}
 Let $B$ be a species of nonseparable graphs and let $C$ denote the species of connected graphs whose blocks are in $B$.
 Then
 \begin{subequations}
   \label{eq:blocks}
   \begin{equation}
     \label{eq:blocksmain}
     B = C \pbrac*{C^{\bullet \abrac{-1}}} + X B' - X
   \end{equation}
   and
   \begin{equation}
     \label{eq:blockssub}
     \specname{E} \pbrac{B'} = \frac{X}{C^{\bullet \abrac{-1}}}.
   \end{equation}
 \end{subequations}
\end{theorem}

It is apparent that the class of nonseparable bipartite graphs is itself exactly the class of blocks that occur in block decompositions of connected bipartite graphs.
We can therefore apply \cref{thm:blocks} to the species $\specname{BP}$ of bipartite blocks.

\begin{theorem}
 \label{thm:znbp}
 Let $\specname{NBP}$ denote the species of $2$-connected bipartite graphs (``bipartite blocks''),  $\specname{CBP}$ the species of connected bipartite graphs, $X$ the species of singletons, and $\con$ the combinatorial logarithm species.
Then $\specname{NBP}$ is determined by
 \begin{subequations}
   \label{eq:nbpexp}
   \begin{equation}
     \label{eq:nbpexpmain}
     \specname{NBP} = \specname{CBP} \pbrac*{\specname{CBP}^{\bullet \abrac{-1}}} + X \cdot \deriv{\specname{NBP}} - X,
   \end{equation}
   where
   \begin{equation}
     \label{eq:nbpexpsub}
     \deriv{\specname{NBP}} = \con \pbrac*{\frac{X}{\specname{CBP}^{\bullet \abrac{-1}}}}.
   \end{equation}
 \end{subequations}
\end{theorem}
We have already calculated the cycle index for the species $\specname{CBP}$, so the calculation of the cycle index of $\specname{NBP}$ is now simply a matter of algebraic expansion.

A generating function for labeled bipartite blocks was given by Harary and Robinson \cite{harrob:bipblocks}, where their analogue of \cref{eq:nbpexp} for the labeled exponential generating function for blocks comes from \cite{forduhl:combprob1}.
However, we could locate no corresponding unlabeled enumeration in the literature.
The numbers of unlabeled nonseparable bipartite graphs for $n \leq 24$ as calculated using our method are given in \cref{tab:bpblocks}, and the Sage code used to compute them is given in \cref{lst:bpcode}.

\appendix

\section{Numerical results}\label{c:enum}
With the tools developed in \cref{c:bpblocks}, we can calculate the cycle indices of the species $\mathcal{NBP}$ of nonseparable bipartite graphs to any finite degree we choose using computational methods.
This result can then be used to enumerate unlabeled bipartite blocks.
We have done so here using Sage \cite{sage} and code listed in \cref{c:code}.
The resulting values appear in \cref{tab:bpblocks}.

\begin{table}[htb]
 \centering
 \caption{Enumerative data for unlabeled bipartite blocks with $n \leq 24$ vertices}
 \label{tab:bpblocks}
 \begin{tabular}{r | r}
   $n$ & Unlabeled bipartite blocks\\\hline
   1 & 1 \\
   2 & 1 \\
   3 & 0 \\
   4 & 1 \\
   5 & 1 \\
   6 & 5 \\
   7 & 8 \\
   8 & 42 \\
   9 & 146 \\
   10 & 956 \\
   11 & 6643 \\
   12 & 65921 \\
   13 & 818448 \\
   14 & 13442572 \\
   15 & 287665498 \\
   16 & 8099980771 \\
   17 & 300760170216 \\
   18 & 14791653463768 \\
   19 & 967055338887805 \\
   20 & 84368806391412395 \\
   21 & 9855854129239183783 \\
   22 & 1546801291978378704267 \\
   23 & 327092325302250220001201 \\
   24 & 93454432085788531687319514
 \end{tabular}
\end{table}
\clearpage

\section{Code listing}\label{c:code}
The functional \cref{eq:nbpexp} characterizes the cycle index of the species $\specname{NBP}$ of bipartite blocks.
In this section we have used the the computer algebra system Sage \cite{sage}
to  adapt the theory into practical algorithms for computing the actual numbers of such structures.
Python/Sage code to compute the coefficients of the ordinary generating function $\widetilde{\specname{NBP}} \pbrac{x}$ of unlabeled bipartite blocks explicitly follows in \cref{lst:bpcode}.

\begin{lstlisting}[caption=Sage code to compute numbers of bipartite blocks (\texttt{bpblocks.sage}), label=lst:bpcode, language=Python, texcl=true]
# ---
# Set up environment
# ---

# Import needed code
from sage.combinat.species.stream import Stream, _integers_from #Infinite generator for lazy power series
from sage.combinat.species.generating_series import CycleIndexSeriesRing
from sage.combinat.species.group_cycle_index_series import GroupCycleIndexSeriesRing
from sage.combinat.species.combinatorial_logarithm import CombinatorialLogarithmSeries

# Set up helper variables
# We'll work with these cycle indices a lot
X = species.SingletonSpecies().cycle_index_series()
E = species.SetSpecies().cycle_index_series()
Omega = CombinatorialLogarithmSeries()

CIS = CycleIndexSeriesRing(QQ) #The ring of cycle index series with rational coefficients
p = SymmetricFunctions(QQ).power() #The ring of symmetric functions (power-sum basis) with rational coefficients

S2 = SymmetricGroup(2) #The group of order 2
GCISR = GroupCycleIndexSeriesRing(S2)
e,t = GCISR.basis().keys()

# Helper method to compute $\specname{X}$ divided by a given cycle index, if possible.
# (Needed for a later calculation where $\specname{F}^{-1}$ is not defined but $\frac{\specname{X}}{\specname{F}}$ is)
def ci_xdiv( f ):
    def p1_dropper( part ):
        assert 1 in part
        return p(part[:-1])
    
    def p1_dropper_sf( sf ):
        assert sf in p
        return p._apply_module_morphism(sf, p1_dropper)
    
    termbuilder = lambda i: CIS([0]*i + [p1_dropper_sf(f.coefficient(i+1)), 0])
    return CIS.sum_generator(termbuilder(i) for i in _integers_from(0))

# ---
# Compute $\gci{\symgp{2}}{\specname{BC}}$
# (This requires a fair amount of code because we need to manually construct the cycle index.)
# ---

# Helper functions for working with partitions
# The union of two partitions is the partition corresponding to their multiset union
def partunion( mu, nu):
    return Partition(sorted(mu.to_list() + nu.to_list(), reverse=true))
    
# For a partition $\mu = \sbrac{\mu_{1}, \mu_{2}, \dots} \vdash m$ and a natural $n$, we define $n \cdot \mu \defeq \sbrac{n \mu_{1}, n \mu_{2}, \dots} \vdash n m$.
# That is, $n \cdot \mu$ is the partition obtained from $\mu$ by multiplying the multiplicity of each part of $\mu$ by $n$.
def partmult( mu, n ):
    return Partition([part * n for part in mu.to_list()])

# Define the cycle index for $\specname{BC}$
# First we compute the number of graphs fixed under $e$ by permutations of cycle
# types $\mu$ and $\nu$ in accordance with \cref{eq:efixedbcgraphs}.
def efixedbcgraphs( mu, nu ):
    return 2**(sum([gcd(i, j) for i in mu for j in nu]))

# Then we build a generator for the terms of the cycle index $\gcielt{\symgp{2}}{\specname{BC}}{e}$.
def egen():
    yield p.zero()
    for n in _integers_from(1):
        yield sum(p(partunion(pair[0], pair[1]))/(pair[0].aut() * pair[1].aut()) * efixedbcgraphs(pair[0], pair[1]) for pair in PartitionTuples(size=n, level=2))

# Then we compute the number of graphs fixed under $\tau$ by a permutation of cycle type
# $\mu$ in accordance with \cref{eq:tfixedbcgraphs}.
def tfixedbcgraphs( mu ):
    return 2**(len(mu) + sum([integer_ceil(m/2) for m in mu]) + sum([gcd(mu[i], mu[j]) for i in range(0, len(mu)) for j in range(i+1, len(mu))]))

# Then we build a generator for the terms of the cycle index $\gcielt{\symgp{2}}{\specname{BC}}{\tau}$.
def tgen():
    yield p.zero()
    for n in _integers_from(1):
        yield p.zero()
        yield sum(tfixedbcgraphs( mu ) * p(partmult(mu, 2))/partmult(mu, 2).aut() for mu in Partitions(n))

# Finally, we use the GroupCycleIndexSeriesRing constructor to define $\gci{\symgp{2}}{\specname{BC}}$.
BC =  GCISR(e)*CIS(egen()) + GCISR(t)*CIS(tgen())

# ---
# Use species algebra to pass through to connected bicolored graphs
# ---

# Define the cycle index for $\specname{CBC}$
CBC = GCISR(Omega).composition(BC)

# Define the cycle indices for $\specname{CBP}$ and $\specname{BP}$ as in \cref{thm:zcbp}.
CBP = CBC.quotient()

# Define the cycle index for $\specname{BP}$ as in \cref{thm:zbp}.
BP = E.composition(CBP)

# ---
# Compute the cycle index for $\specname{NBP}$ using \cref{thm:znbp}
# ---

CBP_pointed_compinv = CBP.pointing().compositional_inverse()
NBP = CBP.composition(CBP_pointed_compinv) + X * Omega.composition((~ci_xdiv(CBP_pointed_compinv))-1)

# Compute the first $k$ coefficients of $\widetilde{\specname{NBP}} \pbrac{x}$
# (Note that this list is zero-indexed!)
print NBP.isotype_generating_series().counts(k)
\end{lstlisting}

\bibliographystyle{amsplain}
\bibliography{sources}

\end{document}